\documentclass[11pt]{article}
\def\bn{\hbox{\it I\hskip -2pt N}}
\def\bz{\hbox{\it Z\hskip -4pt Z}}
\def\bq{\hbox{\it l\hskip -5.5pt Q}}

\def\demo{\noindent{\bf Proof.}}
\newtheorem{theorem}{Theorem}[section]
\newtheorem{lemma}[theorem]{Lemma}
\newtheorem{remark}[theorem]{Remark}
\newtheorem{Corollary}[theorem]{Corollary}
\newtheorem{definition}[theorem]{Definition}

\newtheorem{example}[theorem]{Example}
\begin{document}
\begin{center}
\uppercase{{\bf Complete Intersection Lattice Ideals}}
\end{center}
\advance\baselineskip-3pt
\vspace{2\baselineskip}
\begin{center}
{\sc Marcel Morales and Apostolos Thoma}
\end{center}
\vskip.5truecm\noindent
\vskip.5truecm\noindent
\begin{abstract}
In this paper we completely characterize lattice ideals that
are complete intersections or equivalently complete intersections
finitely generated semigroups of $\bz^n\oplus T$
with no invertible elements, where $T$ is a finite abelian group.  
We also characterize the lattice ideals that are set-theoretic 
complete intersections
on binomials. 
\end{abstract}

\section{Introduction}
 Let $S$ be a  finitely generated, cancellative, abelian semigroup with no invertible elements.
$S$ can be considered as a subsemigroup of a finitely 
generated abelian group $\bz^n\oplus T$ such that $S\cap (-S)=\{{\bf 0}\}$, where $T$ is a torsion group. 
In the case that the torsion group 
is trivial the semigroup $S$ is called {\em affine semigroup}. 
Let $A=\{{\bf a}_i|i\in \{1,\ldots ,m\}\}$ be a set of generators for the semigroup $S$, thus $S=\bn A$,
where $\bn $ is the set of nonnegative integers. Let
$L$ denote the kernel of the group homomorphism from $\bz ^m$ to  $\bz^n\oplus T$ which sends
${\bf e}_i$ to ${\bf a}_i$, where $\{{\bf e}_i|i\in \{1,\ldots ,m\}\}$ is the canonical basis of $\bz^m$.
$L$ is a sublattice of $\bz ^m$, the {\em lattice ideal} associated to $L$ is the binomial ideal
$$I_L= (\{{\bf x}^{{\bf \alpha }^+}- {\bf x}^{{\bf \alpha }^-}|{\bf \alpha } 
={\bf \alpha } ^+-{\bf \alpha} ^- \in L\})\subset K[x_1,\dots ,x_m],$$ where $K$ is a field of any characteristic.
The semigroup $S$ is a complete intersection if and only if $I_L \subset \bq [x_1,\dots ,x_m]$ is a complete intersection, 
which means that the minimal number of generators of $I_L$ is equal to the height of $I_L$.\\
 The problem of determining complete intersection semigroups or equivalently complete
intersection lattice ideals has a long history. It was solved for affine semigroups 
gradually in a series of papers by J. Herzog \cite{H}, Ch. Delorme \cite{Del}, 
R. P. Stanley \cite{Sta}, M. N. Ishida \cite{I}, K. Watanabe \cite{W},  H. Nakajima \cite{Na},
Sch\"afer \cite{S}, J. C. Rosales and P. A. Garcia-Sanchez \cite{R-G}. Finally in 1997  K.G. Fischer, W. Morris
and J. Shapiro \cite{FMS} characterized all complete intersections affine semigroups 
of $\bz^n$ using mixed dominating matrices 
and the notion of semigroup gluing introduced by J. C. Rosales \cite{R}. Recently D. Dais and M. Henk \cite{DH} used Nakajima's classification to describe the precise form
of the binomial equations which determine toric locally complete intersection singularities. \\ 
 Another related problem that drew the attention of a number of authors over the last years 
was the generation of a lattice 
ideal by binomials up to radical \cite{E-S,B-G-N,E-V,BMT,Eto,Eto2,BL}.
In 2002 K. Eto \cite{Eto}
 has characterized complete intersection finitely generated, abelian semigroups with no invertible elements
 or equivalently complete intersection lattice ideals as those that are set-theoretic complete intersection on binomials
in characteristic zero. A generalization of the corresponding result for affine semigroups or
equivalently toric varieties, which was provided by M. Barile, M. Morales and A. Thoma \cite{BMT}. Note that
a binomial ideal $I$ is set-theoretic complete intersection on binomials if there exist $r=height(I)$
binomials $F_1,\dots ,F_r$ such that $rad(I)=rad(F_1,\dots ,F_r)$.
Recently  M. Barile and G. Lyubeznik  \cite{BL} used p-gluing of affine semigroups and {\' e}tale cohomology 
to give a class of toric varieties which are set-theoretic complete intersections only over fields
of one positive characteristic {\em p}.
 
The aim of this article is twofold. On the one hand we give a complete characterization of 
 all finitely generated, cancellative, abelian semigroups with no invertible elements
or equivalently lattice ideals that are complete
intersections by introducing the notion of gluing lattices and extending the notion of semigroup gluing.  
On the other hand we characterize all lattice ideals that are set-theoretic complete intersection on binomials 
in any characteristic by extending the notion of p-gluing. The characterization depends on the characteristic. 

 \section{ Semigroup and lattice gluing }
 \par 
A {\em lattice} is a finitely generated free abelian group.
 A {\em partial character} $(L,\rho )$ on $\bz ^m$ is a 
homomorphism $\rho $ from a sublattice $L$ of 
$\bz^m$ to the multiplicative group $K^*=K-\{ 0 \}$.
 Given a partial character  $(L, \rho )$ on $\bz^m$, we define 
the ideal  
$$I_{L,\rho }:=(\{{\bf x}^{{\bf \alpha }^+}-
\rho({\bf \alpha } ){\bf x}^{{\bf \alpha } ^-}|{\bf \alpha } 
={\bf \alpha } ^+-{\bf \alpha} ^- \in L\})\subset K[x_1,\dots ,x_m]$$
called {\em lattice ideal}.  Here ${\bf \alpha} ^+\in \bn^m$ and  
${\bf \alpha } ^-\in \bn^m$ denote the 
positive and negative part of ${\bf \alpha } $,
respectively, and ${\bf x}^{{\bf \beta }}=x_1^{b_1}\cdots x_m^{b_m}$ 
for ${\bf \beta } =(b_1,\dots ,b_m)\in \bn^m$. We will denote by $F({\bf \alpha })$ the
binomial ${\bf x}^{{\bf \alpha }^+}-{\bf x}^{{\bf \alpha } ^-}$ and by $F_{\rho }({\bf \alpha })$
the binomial ${\bf x}^{{\bf \alpha }^+}-
\rho({\bf \alpha } ){\bf x}^{{\bf \alpha } ^-}$.
Lattice ideals are 
binomial ideals. The 
theory of binomial ideals was developed by Eisenbud and Sturmfels in 
\cite{E-S}. A prime lattice ideal is called a {\em toric ideal}, while the set of 
zeroes in $K^m$ 
is an {\em affine toric variety} in the sense of \cite{St}, since we 
do not require normality. \\
 
Let $A=\{{\bf a}_i| 1\leq i\leq m 
\}\subset \bz^n\oplus T$ be such that the semigroup $\bn A$ has no invertible element.
That means that although the group $\bz^n\oplus T$ has torsion elements,  no nonzero element
in the semigroup $\bn A$ is a torsion element. This remark will be very useful in the sequel.\\
 Let $ \psi : \bz^m \rightarrow \bz^n \oplus T $ be a group homomorphism such that 
$\psi ({\bf e}_i)={\bf a}_i\in \bz^n \oplus T$, where  ${\bf 
e}_1,\dots ,{\bf e}_m$ is the canonical basis of $\bz^m$. We will denote by $L$ the lattice
$ker (\psi )$. The fact that the semigroup $\bn A$ has no invertible element is equivalent
with the fact that the lattice $L$ is {\em positive}, that is $L\cap \bn^m=\{{\bf 0}\}$.
This means that the 
lattice ideal $I_{L,\rho }$ is
homogeneous with respect to some positive grading. In this case by the graded Nakayama's Lemma 
all minimal systems of generators of the
ideal $ I_{L,\rho }$ have the same cardinality.\\

For a lattice $L$ and a prime number $p$, let $(L:p^\infty)$ be the lattice 
$$\{{\bf u}\in \bz^m|p^k{\bf u}\in L \ for \ some \ k\in \bn\}.$$
For a semigroup $S$,   $(S:p^\infty)$ denotes the semigroup 
$$  \{{\bf b}\in \bz^n\oplus T|p^k{\bf b}\in S \ for \ some \ k\in \bn\}.$$
Let   ${ E}\subset \{ 1,\dots ,m \}$, for a set $P\subset \bz$ we denote by 
$$P^E:=\{(p_1, \dots ,p_m)\in \bz^m | p_i\in P\ for \ i\in E, p_i=0\ for \ i\notin { E} \}.$$
  $L_E$  denotes the lattice $L\cap \bz ^{ E} $
and $\bn A^E$ the semigroup generated by $A^E=\{ {\bf a}_i| i\in E\}$.
For a single element ${\bf u}\in \bz^m$ we denote $${\bf u}^E=\{({\bf u}'_i)\in \bz^m|{\bf u}'_i={\bf u}_i
\ for \ i\in E, {\bf u}'_i=0 \ for \ i\notin E\}.$$
\begin{lemma} Let $U\subset L\subset \bz ^m$ be two lattices. Then $p^kL\subset U$ for some $k\in \bn $
if and only if $$(L:p^\infty)=(U:p^\infty).$$
\end{lemma}
\demo Suppose that $p^kL\subset U$ for some $k\in \bn $. From $U\subset L $ we have \\
$(U:p^\infty)\subset (L:p^\infty)$. Let ${\bf u} \in (L:p^\infty)$. Then there exists $n\in \bn$ such that
$p^n{\bf u}\in L$, and the hypothesis implies that $p^{n+k}{\bf u}\in U$. Therefore ${\bf u}\in (U:p^\infty)$.
For the converse, suppose that $L=\sum _{i=1}^r \bz {\bf u}_i$. Then, since ${\bf u}_i\in L$, we have
$${\bf u}_i\in (L:p^\infty)=(U:p^\infty).$$ Which means that there exists $k_i\in \bn $ 
such that $p^{k_i} {\bf u}_i \in U$, $1\leq i\leq r$. By choosing $k$ the maximum of all $k_i$ we have $p^kL\subset U$. 

We give the definitions of semigroup gluing (resp. p-gluing) for subsemigroups of 
 $\bz^n\oplus T$ and gluing (resp. p-gluing) of lattices.

\begin{definition} Let $E_1$, $E_2$ be two nonempty subsets of $\{1, \ldots ,m\}$ such that 
$E_1\cup E_2=\{1, \ldots ,m\}$ and $E_1\cap E_2=\emptyset $. The  semigroup $\bn A$ is called
the gluing (resp. the p-gluing) of the semigroups $\bn A^{E_1}$ and $\bn A^{E_2}$ 
if there is a nonzero 
${\bf a}\in \bn A^{E_1}\cap \bn A^{E_2}$ (resp. ${\bf a}\in 
((\bn A^{E_1}\cap \bn A^{E_2}):p^\infty ))$
such that $\bz{\bf a}=\bz A^{E_1}\cap \bz A^{E_2}$. 
\end{definition}
\begin{definition} Let $E_1$, $E_2$ be two nonempty subsets of $\{1, \ldots ,m\}$ such that 
$E_1\cup E_2=\{1, \ldots ,m\}$ and $E_1\cap E_2=\emptyset $. The lattice $ L$ is called
the gluing (resp. p-gluing) of the lattices $ L_{E_1}$ and $L_{E_2}$ if there is a nonzero ${\bf u}\in L$ 
with ${\bf u}^+={\bf u}^{E_1}$ and ${\bf u}^-=-{\bf u}^{E_2}$, 
such that $L=L_{E_1}+L_{E_2}+<{\bf u}>$ 
(resp. $$(L:p^\infty)=((L_{E_1}+L_{E_2}+<{\bf u}>):p^\infty)).$$
\end{definition}
A set of elements ${\bf a}_1,\dots ,{\bf a}_s$ of $\bz^n\oplus T$ is called linearly
independent if the space of relations is $\{ {\bf 0} \}$, that means the 
relation $\sum _{i=1}^sn_i{\bf a}_i=0$ in $\bz^n\oplus T$, with $n_i\in \bz$, implies $n_1=\cdots =n_s=0$.
\begin{definition}
We call a semigroup  completely glued (resp. p-glued) if it belongs to ${\rm C}$ (resp. P),
which is the smallest class of finitely generated, cancellative, abelian 
semigroups with no invertible elements
that includes all semigroups generated by linearly
independent elements  and is closed under gluing (resp. p-gluing).
\end{definition}
In the sequel we prove some general results that relate the gluing of semigroups with the gluing 
of lattices. We remind the reader that 
$L$ denotes the kernel of the group homomorphism from $\bz ^m$ to  $\bz^n\oplus T$ which sends
${\bf e}_i$ to ${\bf a}_i$, where $\{{\bf e}_i|i\in \{1,\ldots ,m\}\}$ is the canonical basis of $\bz^m$.
Thus with every semigroup $\bn A\subset \bz^n\oplus T$ we associate a lattice $L\subset \bz ^m$. 
Also with every lattice $L\subset \bz ^m$ we associate the semigroup generated by ${\bf e}_i+L$ in $\bz^m/L$,
where $i\in \{1,\dots ,m\}$. We define a lattice to be {\em completely glued} (resp. p-glued) if and only if
the associated semigroup is completely glued (resp. p-glued).

\begin{theorem} The semigroup $\bn A$ is the p-gluing (resp. gluing) of the semigroups $\bn A^{E_1}$ and $\bn A^{E_2}$ 
if and only if the lattice $ L$ is the p-gluing (resp. gluing) of the lattices $ L_{E_1}$ and $ L_{E_2}$. 
\end{theorem}
\demo 
Suppose that $\bn A$ is the p-gluing of $\bn A^{E_1}$ and $\bn A^{E_2}$. Let 
$${\bf a}\in (\bn A^{E_1}\cap \bn A^{E_2}):p^\infty )$$ 
such that $\bz{\bf a}=\bz A^{E_1}\cap \bz A^{E_2}$. Then $p^k{\bf a}=\sum_{i\in E_1}u_i{\bf a}_i=
\sum_{i\in E_2}(-u_i{\bf a}_i)$, for some $k\in \bn$. Then ${\bf u}=(u_i)\in L$ with ${\bf u}^+={\bf u}^{E_1}$ and ${\bf u}^-=-{\bf u}^{E_2}$.
Let ${\bf l}=(l_i)\in (L:p^\infty )$. Then $p^s{\bf l}\in L$ for some $s\in \bn$, which implies 
$\sum_{i\in \{ 1,\dots ,m \}}p^sl_i{\bf a}_i={\bf 0}$. Consider the element
$${\bf b}=\sum_{i\in E_1}p^sl_i{\bf a}_i=\sum_{i\in E_2}(-p^sl_i{\bf a}_i)\in 
\bz A^{E_1}\cap \bz A^{E_2}=\bz{\bf a}.$$
There exists a $\mu \in \bz$ such that ${\bf b}=\mu {\bf a}$, which means 
$\sum_{i\in E_1}p^{k+s}l_i{\bf a}_i=\mu \sum_{i\in E_1}u_i{\bf a}_i$ and 
$\sum_{i\in E_2}(-p^{k+s}l_i{\bf a}_i)=\mu \sum_{i\in E_2}(-u_i{\bf a}_i)$. 
Therefore ${\bf l}_1=(p^{k+s}l_i-\mu u_i)^{E_1}\in L_{E_1}$
and ${\bf l}_2=(p^{k+s}l_i-\mu u_i)^{E_2}\in L_{E_2}$, and $p^{k+s}{\bf l}={\bf l}_1+{\bf l}_2+\mu {\bf u}$.
Therefore  $$(L:p^\infty)\subset((L_{E_1}+L_{E_2}+<{\bf u}>):p^\infty).$$ The other inclusion is obvious. \\
Suppose that $$(L:p^\infty)=((L_{E_1}+L_{E_2}+<{\bf u}>):p^\infty),$$
with ${\bf u}^+={\bf u}^{E_1}$ and ${\bf u}^-=-{\bf u}^{E_2}$. By virtue of Lemma 2.1 there 
exists an $s\in \bn$ such that
$p^sL\subset (L_{E_1}+L_{E_2}+<{\bf u}>)$.
Set ${\bf c}=\sum_{i\in E_1}u_i{\bf a}_i=
\sum_{i\in E_2}-u_i{\bf a}_i$. Then ${\bf c}\in \bn A^{E_1}\cap \bn A^{E_2}$. 
Let ${\bf b}\in \bz A^{E_1}\cap \bz A^{E_2}$, then 
${\bf b}=\sum_{i\in E_1}l_i{\bf a}_i=\sum_{i\in E_2}-l_i{\bf a}_i$.
This implies that ${\bf l}=(l_i)\in L$, therefore $p^s{\bf l}={\bf l}_1+{\bf l}_2+\mu {\bf u}$ for
some ${\bf l}_1\in L_{E_1}$, ${\bf l}_2\in L_{E_2}$ and $\mu \in \bz$. 
But then 
$$p^s{\bf b}=\sum_{i\in E_1}p^sl_i{\bf a}_i=\sum_{i\in E_1}({\bf l}_1+\mu {\bf u}^+)_i{\bf a}_i=\mu {\bf c}. $$
Among the elements of $\bz A^{E_1}\cap \bz A^{E_2}$ choose ${\bf a}$  such that $\mu $ is positive and the smallest possible, set $\mu =\mu _{\bf a}$. 
Then it follows that $\bz{\bf a}=\bz A^{E_1}\cap \bz A^{E_2}$. Now ${\bf c}\in  
\bz A^{E_1}\cap \bz A^{E_2}$, therefore there exists a natural number $\lambda $ such that
${\bf c}=\lambda {\bf a}$. Then from $p^s{\bf a}=\mu _{\bf a}{\bf c}$ we have $p^s{\bf a}=\mu _{\bf a}\lambda {\bf a}$.
Which implies that $\lambda =p^k$ for some $k\in \bn$, since the order of ${\bf a}$ is not finite, as for every nonzero element in $\bn A$. 
Therefore $${\bf a}\in ((\bn A^{E_1}\cap \bn A^{E_2}):p^\infty ).$$
 The proof of the gluing part of the theorem follows from the proof of the p-gluing part by setting $p=1$. 
Actually the second part of the proof is much simpler.
  
The next theorem shows how the  gluing (resp. p-gluing) of lattices reflects on the (resp. radical of the) lattice ideal. 
The first part of the theorem is a generalization of the corresponding result by J. C. Rosales \cite{R}
for toric ideals.
\begin{theorem} Let $E_1$, $E_2$ be two nonempty subsets of $\{1, \ldots ,m\}$ such that 
$E_1\cup E_2=\{1, \ldots ,m\}$ and $E_1\cap E_2=\emptyset $. The lattice 
$L$ is the gluing of the lattices $L_{E_1}$ and $L_{E_2}$ if and only if
$$I_L=I_{L_{E_1}}+I_{L_{E_2}}+<F({\bf u})>,$$ where  ${\bf u}$ is a nonzero element in $L$ such that
 ${\bf u}^+={\bf u}^{E_1}$ and  ${\bf u}^-=-{\bf u}^{E_2}$. 
The lattice 
$L$ is the p-gluing of the lattices $L_{E_1}$ and $L_{E_2}$ if and only if
$$rad(I_L)=rad(I_{L_{E_1}}+I_{L_{E_2}}+<F({\bf u})>),$$ in characteristic $p>0$, 
where  ${\bf u}$ is a nonzero element in $L$ such that
 ${\bf u}^+={\bf u}^{E_1}$ and  ${\bf u}^-=-{\bf u}^{E_2}$. 
\end{theorem}
\demo We prove only the second claim, since the proof of the first is simpler and follows 
from the proof
of the second by putting $p=1$, even in positive characteristic, and taking out the radicals. Suppose that the  lattice 
$L$ is the p-gluing of the lattices $L_{E_1}$ and $L_{E_2}$. Then $$(L:p^\infty)=((L_{E_1}+L_{E_2}+<{\bf u}>):p^\infty).$$
By Theorem 2.5 the semigroup $\bn A$ is the p-gluing of the semigroups 
$\bn A^{E_1}$ and $\bn A^{E_2}$. Then we know that 
$\bz{\bf a}=\bz A^{E_1}\cap \bz A^{E_2}$, where 
$p^k{\bf a}=\sum_{i\in E_1}u_i{\bf a}_i=
\sum_{i\in E_2}(-u_i{\bf a}_i)$.
Let $F({\bf v})\in I_L$. Then ${\bf v}\in L$ and so $\sum_{i=1}^mv_i{\bf a}_i=0$. Then 
$$\sum_{i\in E_1}v_i^+{\bf a}_i+\sum_{i\in E_2}v_i^+{\bf a}_i=\sum_{i\in E_1}v_i^-{\bf a}_i+\sum_{i\in E_2}v_i^-{\bf a}_i.$$
Therefore 
$${\bf \gamma }:=\sum_{i\in E_1}v_i^+{\bf a}_i-\sum_{i\in E_1}v_i^-{\bf a}_i=\sum_{i\in E_2}v_i^-{\bf a}_i-\sum_{i\in E_2}v_i^+{\bf a}_i
\in \bz A^{E_1}\cap \bz A^{E_2}=\bz{\bf a}.$$
That means that ${\bf \gamma }=\tau \sum_{i\in E_1}u_i{\bf a}_i=
\tau \sum_{i\in E_2}(-u_i{\bf a}_i)$, for some $\tau \in \bz$, 
which without loss of generality we can suppose to be
 positive. Then,  since the characteristic is $p>0$, we have $(F({\bf v}))^{p^k}=F(p^k{\bf v})=
{\bf x}^{p^k{\bf v}^+}-{\bf x}^{p^k{\bf v}^-}=
({\bf x}^{p^k{({\bf v}^+)}^{E_1}}-{\bf x}^{p^k{({\bf v}^-)}^{E_1}+\tau {\bf u}^{E_1}}){\bf x}^{p^k{({\bf v}^+)}^{E_2}}
-({\bf x}^{p^k{({\bf v}^-)}^{E_2}}-{\bf x}^{p^k{({\bf v}^+)}^{E_2}+\tau {\bf u}^{E_2}}){\bf x}^{p^k{({\bf v}^-)}^{E_1}}
+{\bf x}^{p^k{({\bf v}^-)}^{E_1}+p^k{({\bf v}^+)}^{E_2}}({\bf x}^{\tau {\bf u}^{E_1}}-{\bf x}^{\tau {\bf u}^{E_2}}).$
From which it is easy to see that $$(F({\bf v}))\in rad(I_{L_{E_1}}+I_{L_{E_2}}+<F({\bf u})>).$$
The reverse inclusion is obvious. \\
Suppose that $rad(I_L)=rad(I_{L_{E_1}}+I_{L_{E_2}}+<F({\bf u})>).$ Let $U$ be the lattice 
$L_{E_1}+L_{E_2}+<{\bf u}>$ then $U\subset L$ and thus also $I_U\subset I_L$ .
Also note that $I_{L_{E_1}}\subset I_U$, $I_{L_{E_2}}\subset I_U$ and $<F({\bf u})>\subset I_U$. Therefore
$I_{L_{E_1}}+I_{L_{E_2}}+<F({\bf u})>\subset I_U$, which implies 
$$rad(I_{L_{E_1}}+I_{L_{E_2}}+<F({\bf u})>)\subset rad(I_U).$$ Then from the hypothesis we have  $rad(I_U)=rad(I_L)$.
It follows from \cite{E-S}, Corollary  2.2, that in characteristic zero $I_U=I_L$ and so $U=L$, and in characteristic
$p>0$ that $I_{(U:p^\infty )}=I_{(L:p^\infty )}$ and so $(U:p^\infty )=(L:p^\infty )$. Note that in \cite{E-S} 
$(L:p^\infty )$ is denoted by $Sat_p(L)$.

\section{ Complete Intersections }
In this section  we will give a series of results that will characterize complete intersection lattice ideals
 and complete intersection semigroups. We also characterize lattice ideals that are set-theoretic complete intersections
on binomials.\\
Let $L$ be a nonzero positive sublattice of $\bz^m $ of 
rank $ r $, and $(L, \rho )$ 
be a partial character on $\bz^m $.  The height of the lattice ideal $I_{L, \rho}$ is equal
to $r$, the rank of the lattice $L$, see \cite{E-S}, Corollary 2.2.

\begin{remark}  Any variable $x_{i}$ is a nonzero divisor for $I_{L,\rho }$.

\end{remark} 
We grade $K[x_1,\cdots ,x_m]$ by setting $deg_{(\bz ^m/L)}(x_i)={\bf a}_i$,
for $i\in \{1, \ldots ,m\}$. Then the $\bz ^m/L$-degree of the monomial $x^{\bf u}$ is
$$deg_{(\bz ^m/L)}(x^{\bf u})=u_1{\bf a}_1+\cdots +u_m{\bf a}_m\in \bn A,$$
where $\bn A$ is the semigroup generated by $A$. The lattice ideal 
$I_{L,\rho }\subset K[x_{1},\ldots,x_{m}]$ is $\bz ^m/L$-homogeneous, since all generators
are $\bz ^m/L$-homogeneous. In particular, let $ {\bf v}\in \bz^m, A,B \in K^* $ and $G({\bf v})=A{\bf 
x}^{{\bf v}^+}- B{\bf x}^{{\bf v} ^-}$, then $G({\bf v})\in I_{L,\rho }$ 
implies ${\bf v}\in L$. Since, if  ${\bf v}\notin L$, then $G({\bf v})$ is not 
$\bz ^m/L$-homogeneous.  Then  the monomial ${\bf 
x}^{{\bf v}_+}$ must be in $ I_{L,\rho }$ since $I_{L,\rho }$ is 
$\bz ^m/L$-homogeneous.  This is impossible since any variable $x_{i}$ is a 
nonzero divisor for $ I_{L,\rho }$.

 \begin{lemma} Let $I,J,K\subset R$ be three ideals in a 
noetherian ring $R$ such that $J\subset I$ and $rad(I)=rad (J),$
then $$rad(I+K)=rad (J+K).$$
\end{lemma}
 \demo The inclusion $rad (J+K)\subset rad(I+K)$ is clear. Now 
let $ g\in rad(I+K)$. Then $g^q\in I+K$ and we can write 
$g^q=h_{1}+h_{2}$, with $h_{1}\in I, h_{2}\in K$. Hence there exists $l$ 
such that $h_{1}^l\in J$, so $g^{ ql }=h_{1}^l+h'_{2}$ with $h'_{2}\in 
K$, which proves the assertion.

\begin{lemma} Consider $ r $ vectors ${\bf u}_1, 
\ldots, {\bf u}_r\in \bz^m $, let $L=\sum _{i=1}^r\bz u_i$ be the lattice generated by them.
The following are equivalent:
 \begin{enumerate}
 	\item  $I_L=(F({\bf u}_1),\ldots,F({\bf u}_r))$ and 
$F({\bf u}_1),\ldots,F({\bf u}_r)$ is a 
regular sequence,
 
 \item  $I_{L,\rho}=(F_{\rho }({\bf u}_1),\ldots,F_{\rho }({\bf u}_r))$ and 
 $F_{\rho }({\bf u}_1),\ldots,F_{\rho }({\bf u}_r)$ is a regular sequence for any 
 partial character $(L, \rho )$ on $\bz^m $. 
 \end{enumerate}
\end{lemma}
\demo  First we remark that any variable 
 $x_{i}$ is a non zero divisor of $I_{L}$, 
 this implies that the sequence $F({\bf u}_1),\ldots,F({\bf u}_r),x_{1}\ldots x_{m}$, is a 
 regular sequence. Let $(L, \rho )$ be a partial character on $\bz^m $.
Then $\rho ({\bf u})$ is a unit for every ${\bf u}\in L$. Thus by \cite{SSS}, Theorem 2.7, the sequence 
$F_{\rho }({\bf u}_1),\ldots,F_{\rho }({\bf u}_r),x_{1}\ldots 
 x_{m},$ is regular.
Let ${\bf u}\in L$ any nonzero vector, we can write 
${\bf u}=n_{1}{\bf u}_{1}+ \ldots+n_{r}{\bf u}_{r}$. From the identity
$$\frac{{{\bf x}^{{\bf u}^+}}}{{{\bf x}^{{\bf u}^-}}}-\rho ({\bf u})=\prod _{i=1}^r
(\frac{{{\bf x}^{{\bf u}_i^+}}}{{{\bf x}^{{\bf u}_i^-}}})^{n_i}-\prod _{i=1}^r 
(\frac{{\rho ({\bf u}_i^+})}{{\rho ({\bf u}_i^-})})^{n_i}$$
by clearing denominators we get an identity in $K[x_{1},\ldots ,x_{m}]$ which shows that 
there exists  a monomial $P$ such that 
  $PF_{\rho }({\bf u})$ belongs to 
$(F_{\rho }({\bf u}_1),\ldots,F_{\rho }({\bf u}_r))$. But 
 $F_{\rho }({\bf u}_1),\ldots,F_{\rho }({\bf u}_r),x_{1}\ldots x_{m}$, is a regular sequence which 
 implies that $F_{\rho }({\bf u})\in (F_{\rho }({\bf u}_1),\ldots,F_{\rho }({\bf u}_r))$,
therefore $I_{L,\rho}=(F_{\rho }({\bf u}_1),\ldots,F_{\rho }({\bf u}_r))$.
\\ The 
 proof of the other implication follows from applying $(2)$ to the trivial character.

 \begin{Corollary} For any lattice ideal $I_{L,\rho}$ the fact that 
 $I_{L,\rho}$ is a complete intersection is independent from the 
 character $\rho$.
 	
 \end{Corollary}
 \begin{definition} \cite{FS} A matrix $ M $ with coefficients in $ \bz 
$ is called mixed if every row has a positive and a negative entry.  $ 
M $ is called dominating if it does not contain any square mixed 
submatrix.\\
We also define the empty matrix ($0\times d$) to be mixed dominating.

\end{definition}

We denote by $M({\bf u}_{1}, \ldots, {\bf u}_{r})$ the $r\times m$ matrix whose rows are the vectors  
${\bf u}_{1}, \ldots, {\bf u}_{r}$ of $\bz ^m$. 
 
  \begin{theorem} Let $L$ be a nonzero positive sublattice of $\bz^m$ of rank $ r $, 
and $(L, \rho )$ be a partial character on $\bz^m$.
Consider $ r $ vectors $ {\bf u}_{1}, \ldots, {\bf u}_{r}\in L $.
The 
following are equivalent:

 \begin{enumerate}
 \item $ rad(I_{L,\rho}) = rad(F_{\rho}({\bf u}_1),\ldots,F_{\rho}({\bf u}_r)) $.
 \item 
 	\begin{itemize} 	
 		\item  the matrix $M({\bf u}_{1}, \ldots, {\bf u}_{r})$ is mixed dominating,
 	
 		\item  in characteristic $0$ we have that $L= \sum _{i=1}^r \bz {\bf u}_i$ and in characteristic $p>0$,  
$$(L:p^\infty)=(\sum _{i=1}^r \bz {\bf u}_i:p^\infty).$$

 	\end{itemize} 
\end{enumerate}
\end{theorem}

\demo ({\bf 1}  $\Rightarrow $ {\bf 2}) 
 Since  $L \subset \bz^m $ is a positive 
sublattice, the matrix $M$ is mixed. 
 Now we prove that $ M $ is dominating, i.e., no square 
submatrix of $ M $ is mixed.    Assume that $ N $ is a mixed $ s\times s $ submatrix of 
$ M $, with $ s\geq 1 $ and suppose that $ s $ is maximal with respect 
to this property.  Then up to permutations of the rows and of the 
variables we may assume that $ N $ consists of the first $ s $ lines 
and the first $ s $ columns, so that we can write: $$ M = \pmatrix{ N 
& B\cr C & D\cr }.$$
From Lemma 3.2 we have  $$rad(I_{L}+(x_{1},\ldots,x_{s}))=rad 
(F({\bf u}_1),\ldots,F({\bf u}_r),x_{1},\ldots,x_{s}).$$ 

Since $N$ is mixed, $$(F({\bf u}_1),\ldots,F({\bf u}_s))\subset 
(x_{1},\ldots,x_{s})$$ so in fact we have 
$$rad(I_{L}+(x_{1},\ldots,x_{s}))=rad 
(F({\bf u}_{s+1}),\ldots,F({\bf u}_r),x_{1},\ldots,x_{s}).$$ On the other 
hand $x_{1}$ is not a zero divisor of $I_{L}$, therefore $height 
(rad(I_{L}+(x_{1},\ldots,x_{s}))) \geq r+1 $, but the height of $rad 
(F({\bf u}_{s+1}),\ldots,F({\bf u}_r), x_{1},\ldots,x_{s})$ is atmost $r$.  
This is a contradiction, therefore $M$ is mixed dominating.
 
Since $M({\bf u}_1,\ldots,{\bf u}_r)$ is mixed dominating,
 by Fischer-Shapiro \cite{FS}, Theorem  2.9, we get that the ideal
 $(F({\bf u}_1),\ldots,F({\bf u}_r))$ is equal to the lattice ideal $I_U$, 
where $U=\sum _{i=1}^{r} \bz {\bf u}_i$. By Lemma 3.3 this implies
  $(F_{\rho}({\bf u}_1),\ldots,F_{\rho}({\bf u}_r))=I_{U,\rho}$.
Now by hypothesis there exists 
  $k$ such that $$F_{\rho}({\bf v})^{ p^k }\in (F_{\rho}({\bf u}_1),\ldots,F_{\rho}({\bf u}_r))$$ for any ${\bf v}\in 
  L$. $(L, \rho )$ 
is  a partial character on $\bz^m $ therefore $\rho(p^k{\bf v})=(\rho({\bf v}))^{p^k}$. 
If the characteristic of $K$ is equal to $p$, 
this implies $F_{\rho}({\bf v})^{ p^k }=F_{\rho}(p^k{\bf v})\in 
  (F_{\rho}({\bf u}_1),\ldots,F_{\rho}({\bf u}_r))=I_{U,{\rho}}$ and then $p^k{\bf v} \in U$,
 since $I_{U,{\rho}}$ is
$\bz ^m/U$-homogeneous.  Therefore $(L:p^\infty)=(U:p^\infty)$.
If the characteristic of $K$ is zero, $I_{U,{\rho}}$ is a radical ideal, see Eisenbud-Sturmfels \cite{E-S}, Corollary 2.2,
then $$F_{\rho}({\bf v})^{ 
  p^k }\in (F_{\rho}({\bf u}_1),\ldots,F_{\rho}({\bf u}_r))=I_{U,{\rho}}$$ implies $F_{\rho}({\bf v})\in 
  I_{U,{\rho}}$, therefore ${\bf v} \in U$ and $L=U$.

({\bf 2}  $\Rightarrow $ {\bf 1}) Since $M({\bf u}_1,\ldots,{\bf u}_r)$ is mixed dominating,
 by Fischer-Shapiro \cite{FS}, Theorem  2.9, we get 
 $(F({\bf u}_1),\ldots,F({\bf u}_r))=I_U$ and by Lemma 3.3 this implies
  $(F_{\rho}({\bf u}_1),\ldots,F_{\rho}({\bf u}_r))=I_{U,\rho}$.  If the characteristic of $K$ is zero we have $U=L$, so  
  $I_{L,\rho}= I_{U,\rho}=(F_{\rho}({\bf u}_1),\ldots,F_{\rho}({\bf u}_r))$.  If 
the characteristic 
  of $K$ is $p$ positive, for any ${\bf v}\in L$, we have $$F_{\rho}({\bf v})^{ p^k }=F_{\rho}(p^k{\bf v})\in 
  I_{U,\rho}=(F({\bf u}_1)_{\rho},\ldots,F({\bf u}_r)_{\rho})$$ and then 
  $I_{L,\rho}\subset rad(I_{U,\rho})= 
  rad(F_{\rho}({\bf u}_1),\ldots,F_{\rho}({\bf u}_r))$.  This completes the proof. 

\begin{remark}
By  Fischer-Shapiro \cite{FS}, Corollary 2.8,  if the matrix  
$M({\bf u}_{1}, \ldots, {\bf u}_{r})$ is mixed dominating then the vectors 
${\bf u}_{1}, \ldots, {\bf u}_{r}$ are linearly independent.
\end{remark}
 \begin{Corollary} For any lattice ideal $I_{L,\rho}$ the 
fact that $I_{L,\rho}$ is a set-theoretical complete intersection on binomials is 
independent from the character $\rho$. Moreover, if $ rad(I_L) = rad(F({\bf u}_1),\ldots,F({\bf u}_r))$,
then for any character $\rho$
$$ rad(I_{L,\rho}) = rad(F_{\rho}({\bf u}_1),\ldots,F_{\rho}({\bf u}_r)).$$
 \end{Corollary}
The proof follows from Theorem 3.6, since condition ($2$) is independent of the character.
\begin{theorem}  Let $L$ be a nonzero positive sublattice of $\bz^m 
$ of rank $ r $, and $(L, \rho )$ be a partial character on $\bz^m 
$.  Consider $ r $ vectors $ {\bf u}_{1}, \ldots, {\bf u}_{r}\in L $, 
the 
following are equivalent:

 \begin{enumerate}
 \item $ I_{L,\rho} =(F_{\rho}({\bf u}_1),\ldots,F_{\rho}({\bf u}_r)) $,
 \item 
 	\begin{itemize}
 \item The matrix $M({\bf u}_1,\ldots,{\bf u}_r)$ is mixed dominating,
 	
 		\item   $L= \sum _{i=1}^{r}\bz{\bf u}_i$.

 	\end{itemize}
 	 
 \end{enumerate}
 
 \end{theorem}
The proof follows from the proof of Theorem 3.6 by taking out the radicals and putting $p=1$ even in positive 
characteristic. Theorem 3.9 characterizes complete intersection lattice ideals: a lattice 
ideal $I_{L,\rho}$ is a complete intersection if and only if the lattice $L$ has a basis 
${\bf u}_{1}, \ldots, {\bf u}_{r}$ such that the matrix $M({\bf u}_1,\ldots,{\bf u}_r)$
is mixed dominating.
\begin{Corollary} Let $L$ be a nonzero positive sublattice of $\bz^m 
$ of rank $ r $, and $(L, \rho )$ be a partial character on $\bz^m 
$. If the characteristic of $K$ is zero,  we have $rad(I_{L,\rho }) = rad(F_{\rho}({\bf u}_1),\ldots,F_{\rho}({\bf u}_r))$ if and only if $I_{L,\rho } = (F_{\rho}({\bf u}_1),\ldots,F_{\rho}({\bf u}_r))$. 
 \end{Corollary}
The proof of the Corollary follow from the proof of Theorem 3.6. 
Corollary 3.10 states that in zero characteristic a lattice ideal is a set-theoretic complete intersection on binomials if and only if it is a
complete intersection, see also \cite{Eto}, Theorem 2.1.

The aim of the next theorems is to prove Theorems 3.15 and 3.16, which give an exact characterization of complete
intersection lattice ideals and complete intersection semigroups. Lattices that correspond to lattice ideals that
are set-theoretic complete intersection on binomials are also characterized.
\par We recall the following decomposition theorem of K. Fischer, W. Morris and J. Shapiro, for mixed dominating
matrices (see \cite{FMS}, Theorem 2.2) whose claim we adjust to our notation.

\begin{theorem} Let $M({\bf u}_1,\ldots,{\bf u}_r)$ be a mixed dominating $ r\times m $
matrix with $m\geq r>0.$ Then there exist $E_1$, $E_2$ disjoint nonempty subsets of $\{1, \ldots ,m\}$ with 
$E_1\cup E_2=\{1, \ldots ,m\}$, 
and disjoint subsets $S_1$, $S_2$ of $\{1, \ldots ,r\}$ with 
$S_1\cup S_2=\{1, \ldots ,r\}-\{ q\}$ for some $q$, such that
 the matrices $M(\{{\bf u}_i|i\in S_1\})$, $M(\{{\bf u}_i|i\in S_2\})$ are mixed dominating,
where $({\bf u}_i)^{E_j}={\bf u}_i$ for every $i\in S_j$, $j\in \{1,2\}$ and 
 $({\bf u}_q)^{E_1}={\bf u}_q^+$,  $({\bf u}_q)^{E_2}=-{\bf u}_q^-$. 
\end{theorem}

\begin{lemma} The notation being that of Theorem 3.11 we have for $j\in \{ 1,2\}$,
 $$(\sum_{i=1}^r\bz {\bf u}_i)_{E_j}=\sum_{i\in S_j}\bz {\bf u}_i$$
and the lattice $U=\sum_{i=1}^r\bz {\bf u}_i$ is the gluing of the lattices $U_{E_1}, U_{E_2}$.
\end{lemma}
\demo Without loss of generality we take $j=1$. Recall that $L_{E_1}=L\cap \bz^{E_1}$, and since 
$({\bf u}_i)^{E_1}={\bf u}_i$ for every $i\in S_1$ we conclude that $\sum_{i\in S_1}\bz {\bf u}_i 
\subset (\sum _{i=1}^{r}\bz{\bf u}_i)_{E_1}$.
Let ${\bf u}\in (\sum _{i=1}^{r}\bz{\bf u}_i)_{E_1}\subset 
\sum _{i=1}^{r}\bz{\bf u}_i$. Then ${\bf u}={\bf u}^{E_1}$ and
${\bf u}=\sum_{i\in S_1}\lambda _i {\bf u}_i+\sum_{i\in S_2}
\lambda _i {\bf u}_i+\lambda _q{\bf u}_q$.
From which we have that 
${\bf u}^{E_1}=\sum_{i\in S_1}\lambda _i {\bf u}_i^{E_1}+\sum_{i\in S_2}\lambda _i {\bf u}_i^{E_1}+
\lambda _q{\bf u}_q^{E_1}$. But then 
${\bf u}=\sum_{i\in S_1}\lambda _i {\bf u}_i+\lambda _q{\bf u}_q^+$. The last equality implies that
the vector $\lambda _q{\bf u}_q^+$ belongs to the positive lattice $\sum_{i=1}^r\bz {\bf u}_i$, which is impossible
except if  $\lambda _q=0$. Thus ${\bf u}=\sum_{i\in S_1}\lambda _i {\bf u}_i$.  
We conclude that $(\sum _{i=1}^{r}\bz{\bf u}_i)_{E_1}=\sum_{i\in S_1}\bz {\bf u}_i.$ 
Therefore, $U=U_{E_1}+U_{E_2}+<{\bf u}_q>$, where $({\bf u}_q)^{E_1}={\bf u}_q^+$,  $({\bf u}_q)^{E_2}=-{\bf u}_q^-$.

\begin{theorem} Let $K$ be a field of positive characteristic $p$. The lattice ideal $I_{L,\rho }\subset K[x_1,\dots ,x_m]$
is set-theoretic complete intersection on binomials if and only if the lattice $L$ is the p-gluing of the two
lattices $L_{E_1}$ and $L_{E_2}$ and both lattice ideals $I_{L_{E_1},\rho }, I_{L_{E_2},\rho }$ are set-theoretic complete
intersections on binomials. 
\end{theorem}
\demo
 Suppose that 
$rad(I_{L,\rho })=rad(F({\bf u}_1),\ldots,F({\bf u}_r)))$. Then Theorem 3.6 gives us that
the matrix $M({\bf u}_1,\ldots,{\bf u}_r)$ is mixed dominating.
Therefore there exist $E_1$, $E_2$,  $S_1$, $S_2$ as provided in Theorem 3.11.
By virtue of Lemma 3.12 the lattice $U=\sum_{i=1}^r\bz {\bf u}_i$ is the gluing of the lattices $U_{E_1}, U_{E_2}$.
Now $U\subset L$ and  from Theorem 3.6 we have $$(L:p^\infty)=(\sum _{i=1}^r \bz {\bf u}_i:p^\infty),$$ therefore
by Lemma 2.1 there exists a positive integer $k$ such that $p^kL\subset U=U_{E_1}+U_{E_2}+<{\bf u}_q>$.
But $U_{E_1}\subset L_{E_1}$ and $U_{E_2}\subset L_{E_2}$ so that $p^kL\subset L_{E_1}+L_{E_2}+<{\bf u}_q>$.
Note also that $L_{E_1}+L_{E_2}+<{\bf u}_q>\subset L$. Therefore, by Lemma 2.1, we have that 
$$(L:p^\infty)=(L_{E_1}+L_{E_2}+<{\bf u}_q>:p^\infty).$$ Which means that $L$ is the p-gluing of 
 $L_{E_1}$ and $L_{E_2}$.\\
Note also that $(L_{E_j}:p^\infty)=(L:p^\infty)_{E_j}=(U:p^\infty)_{E_j}=(U_{E_j}:p^\infty)$, for $j\in \{1,2\}$.
By Remark 3.7 the vectors $ {\bf u}_{1}, \ldots, {\bf u}_{r} $ are linearly independent and by 
Theorem 3.11 the matrices $M(\{{\bf u}_i|i\in S_j\})$ are mixed dominating, for $j\in \{1,2\}$.
Therefore, by Theorem 3.6 again, we conclude that  $rad(I_{L_{E_j}})  = rad(F({\bf u}_i)|i\in S_j)$, for $j\in \{1,2\}$.
Recall that $height(I_{L_{E_j},\rho })=rank(\sum_{i\in S_j}\bz {\bf u}_i)=|S_j|$, for $j\in \{1,2\}$, therefore
both $I_{L_{E_1},\rho }, I_{L_{E_2},\rho }$ are set-theoretic complete
intersections on binomials. \\
The proof of the converse implication follows from Theorem 2.6 and the remark that for the lattice (p-) gluing 
for positive lattices we have $rank(L)=rank(L_{E_1})+rank(L_{E_2})+1$.

 Notice that by, Corollary 3.10, in the zero characteristic case, lattice ideals that are binomial set theoretic 
complete intersection  are complete intersections.  Therefore they are characterized also by the next Theorem.
 \begin{theorem} 
The lattice ideal $I_{L,\rho }\subset K[x_1,\dots ,x_m]$
is a complete intersection if and only if the lattice $L$ is the gluing of the two
lattices $L_{E_1}$ and $L_{E_2}$ and both lattice ideals $I_{L_{E_1},\rho }, I_{L_{E_2},\rho }$ are  complete
intersections. 
\end{theorem}
\demo  The proof follows the lines of the proof of Theorem 3.13 by taking out the radicals and putting $p=1$ even in positive 
characteristic.  

The next theorem is the main result of the article and characterizes all lattice ideals that are  complete intersections  and
also all lattice ideals that are set-theoretic complete intersections on binomials, in all
characteristics. 
\begin{theorem} Let $K$ be a field of any characteristic. 
The lattice ideal $I_{L,\rho }\subset K[x_1,\dots ,x_m]$
is a complete intersection if and only if the lattice $L$ is completely glued.\\
In the characteristic zero case (resp. positive characteristic p case), 
the lattice ideal $I_{L,\rho }$
is a set-theoretic complete intersection on binomials if and only if the lattice $L$ is completely glued (resp. completely
p-glued). 
\end{theorem}
\demo  The proof follows by induction on the rank $r$ and is based on Theorems 3.13, 3.14. Note that if a lattice has rank
zero then the elements of the associated semigroup are linearly independent and therefore the lattice is completely (resp. p-) glued and of course a complete intersection. 

The property for a lattice ideal to be a complete intersection does not depend on the field, but only on the lattice $L\subset \bz^m$. Therefore,
translating Theorem 3.15 for semigroups we have:
\begin{theorem} A finitely generated, cancellative, abelian semigroup with no invertible elements is a
 complete intersection if and only if it is completely glued. \end{theorem}
Theorem 3.16 restricted to affine semigroups gives an exact characterization
 of complete intersection affine semigroups: {\em an affine
semigroup  is a complete intersection if and only if it is completely glued}. An affine
semigroup is {\em completely glued} if it belongs to the smallest class of affine semigroups
that includes all free affine semigroups and is closed under gluing.
\begin{example}{\rm The results of this section help us to provide examples of lattice ideals
that are complete intersections or set-theoretic complete intersections on binomials.
Any mixed dominating integer matrix $M({\bf u}_1,\ldots,{\bf u}_r)$ gives 
a completely glued lattice, the  $L=\sum_{i=1}^r\bz {\bf u}_i$, and a complete intersection
lattice ideal, the $I_{L,\rho }$ in $K[x_1,\dots ,x_m]$, where $K$ is any field and $(L,\rho)$
a partial character 
 on $\bz^m$ . Also the semigroup
$<{\bf e}_i+L|i\in \{1,\dots ,m\}>\subset\bz^m/L$ is completely glued. Considering a lattice $L'$ such that 
$(L':p^\infty)=(L:p^\infty)$ for some prime number $p$, the lattice ideal 
$I_{L',\rho }$ in $K[x_1,\dots ,x_m]$ is set-theoretic complete intersection 
on binomials, where $K$ is  a field of characteristic $p$.
\\
Mixed dominating matrices  can be constructed easily. 
Let $M_1$ and $M_2$ be mixed dominating matrices of sizes $m_1\times n_1$ and $m_2\times n_2$ with $m_1\geq 0$
and $m_2\geq 0$. Let ${\bf u^+}\in \bn^{n_1}$ and ${\bf u^-}\in \bn^{n_2}$ be any two vectors. Then the matrix 
$$  \pmatrix{ M_1
& 0 \cr 0 & M_2\cr {\bf u}^+ & {\bf -u}^- \cr},$$
is mixed dominating. To start with, we can consider both  matrices $M_1$, $M_2$ to be empty. Subsequently we use 
already constructed mixed dominating matrices to construct new ones. 
Actually the decomposition theorem, see Theorem 3.11 or \cite{FMS}, Theorem 2.2, of mixed dominating 
matrices says that all mixed dominating matrices can be taken 
in this way.\\
For example, take as $M_1$ the $1\times 3$ mixed dominating matrix 
$  \pmatrix{ 1 & 3 & -4 \cr}$, $M_2$ the empty $0\times 1$ matrix, ${\bf u^+}=(3,1,0)$
and ${\bf u^-}=(4)$. Then the matrix $$  \pmatrix{ 1 & 3 & -4 & 0 \cr
3 & 1 & 0 & -4 \cr}$$ is mixed dominating. Therefore the lattice $L=\bz (1,3,-4,0)+
\bz (3,1,0,-4)$ is completely glued and the lattice ideal $I_{L,\rho }$ is a complete
intersection for any character $\rho $. The associated semigroup of the lattice $L$ is 
isomorphic to the semigroup generated by $(4,0,0)$, $(0,4,0)$, $(1,3,0)$ and $(3,1,1)$
in $\bz^2\oplus \bz _4$. Which is completely glued. \\
Let $L'= \bz (1,3,-4,0)+
\bz (0,2,-3,1)$. $L'$ is the associated lattice of the affine semigroup generated 
by $(4,0)$, $(0,4)$, $(1,3)$ and $(3,1)$ in $\bz^2$. Which is not a complete intersection affine semigroup. Therefore there is no basis ${\bf u}_1, {\bf u}_2$ of $L'$ such that
the matrix $M({\bf u}_1, {\bf u}_2)$ is mixed dominating. Notice that 
$(L':2^\infty)=(L:2^\infty)$, since $4L'\subset L\subset L'$. This implies that in characteristic 2 the two ideals 
$I_{L',\rho }$, $I_{L,\rho }$ have the same radical. Therefore $I_{L',\rho }$ is set-theoretic
complete intersection on binomials in characteristic 2.}
\end{example}

\section{Extreme rays of a complete intersection semigroup cone}

Let $\phi $ be the projection homomorphism from $\bz^n\oplus T$ to $\bz^n$ and denote
$\phi ({\bf b})={\bf \overline  b}$ for ${\bf b}\in \bz^n\oplus T$.
Let $\overline A=\{{\bf \overline  a}_i| 1\leq i\leq m \}$. We associate with the semigroup ${\bn }A$
(or with the lattice 
ideal $I_{L,\rho }$) 
 the rational polyhedral cone  
$\sigma =pos_{\bq}(\overline A):=\{l_1{\bf \overline  a}_1+\cdots +l_m{\bf \overline  a}_m|
l_i\in \bq \ and \ l_i\geq 0\}$.  A cone $\sigma $ is {\em strongly 
convex} 
if $\sigma \cap -\sigma =\{ {\bf 0} \}$. 
The condition that the lattice $L$ is positive is equivalent with 
the condition 
that the cone $\sigma $
is strongly convex.\\
A ray $R$ in the cone of $\overline A$ is an {\em extreme ray} of the cone of $\overline A$, if
given any vector ${\bf u}\in R$, positive integers $\mu ,c_1, \dots ,c_t$ and elements
${\bf w}_1, \dots ,{\bf w}_t$ of ${\bn \overline A}$ such that $$\mu {\bf u}=c_1{\bf w}_1+\dots c_t{\bf w}_t,$$
then ${\bf w}_j\in R$ for all $j=1, \dots ,t$. In \cite{FMS} it was shown that for  an 
$n$-dimensional  complete intersection affine semigroup  with $n\ge 2$, its  cone contains
no more than $2n-2$ extreme rays. The corresponding statement is true for semigroups
of $\bz^n\oplus T$ or equivalently lattice ideals which are complete intersections. But also
for lattice ideals that are set theoretic complete intersections on binomials.
\begin{theorem} Let ${\bn }A$ be an n-dimensional semigroup of $\bz^n\oplus T$ which is completely glued
or completely p-glued, $n\geq 2$.  Then the cone of $\overline A$
contains no more than $2n-2$ extreme rays.
\end{theorem}
\demo The proof almost follows the lines of the proof of Corollary 2.4 in \cite{FMS}. Let ${\bn }A$ be a semigroup of $\bz^n\oplus T$ which is completely glued
or completely p-glued. Let $ \psi : \bz^m \rightarrow \bz^n \oplus T $ be the group homomorphism such that 
$\psi ({\bf e}_i)={\bf a}_i\in \bz^n \oplus T$, where  ${\bf 
e}_1,\dots ,{\bf e}_m$ is the canonical basis of $\bz^m$. Let $L$ be the lattice
$ker (\psi )$ of rank $r=m-n$. We will use induction on $r$. If $r=0$ then $m=n$. Hence
the vectors in ${\overline A}$ are linearly independent and the cone has exactly $n$ extreme rays.
Since $n\geq 2$, we have $n\leq 2n-2$.\\
If $r\geq 1$, we can write $A$ as the disjoint union of $A^{E_1}$, $A^{E_2}$ such that
 $\bz{\bf a}=\bz A^{E_1}\cap \bz A^{E_2}$ and 
there is a multiple of ${\bf a}$ in $\bn A^{E_1}\cap \bn A^{E_2}$, for some  disjoint subsets $E_1, E_2$ of
$\{1, \dots ,m\}$. Then we have ${\bf \overline  a}\in \bz{\overline A}^{E_1}\cap \bz{\overline A}^{E_2}$.
Let ${\bf \overline  b}\in \bz{\overline A}^{E_1}\cap \bz{\overline A}^{E_2}$. Then
$g{\bf b}\in \bz A^{E_1}\cap \bz A^{E_2}=\bz{\bf a}$, where $g$ is the order of the finite
group $T$. Therefore $g{\bf b}=\lambda {\bf a}$ and so $g{\bf \overline b}=\lambda {\bf \overline a}$.
Thus $\bz{\overline A}^{E_1}\cap \bz{\overline A}^{E_2}$ is one dimensional and if
${\bf \overline c}$ is any generator, then ${\bf \overline a}=\mu {\bf \overline c}$. We conclude that
a multiple of ${\bf \overline c}$ belongs to $\bn{\overline A}^{E_1}\cap \bn{\overline A}^{E_2}$.\\
Let $n_1,n_2$ be the dimensions of ${\bn }A^{E_1},{\bn }A^{E_2}$, respectively.
Then $n_1+n_2=n+1$. Let $r_i$ be the rank of the lattice $L_{E_i}$, $i\in \{1,2\}$. 
It follows from $n_1+n_2=n+1$ that $r_1+r_2=r-1$. Therefore each $r_i$ is less than $r$.
Each extreme ray of the cone of ${\overline A}$ is an extreme ray for either the cone of
${\overline A}^{E_1}$ or ${\overline A}^{E_2}$. Therefore, the number of extreme rays
of the cone of ${\overline A}$ is bounded by the sum of the number of extreme rays in the cones
of ${\overline A}^{E_1}$ and ${\overline A}^{E_2}$. Hence as long as $n_i\geq 2$, the inductive
hypothesis gives that the number of extreme rays of the cone of ${\overline A}$
is bounded by $2n_1-2+2n_2-2=2n-2$. But if $r_1=1$ say, then since
the two cones of ${\overline A}^{E_1}$ and ${\overline A}^{E_2}$ intersect in a semiline, it
follows that the cone of ${\overline A}^{E_1}$  is contained in the cone of 
${\overline A}^{E_2}$. Therefore the cone of ${\overline A}$ is the same with the
cone of ${\overline A}^{E_2}$. But $r_2$ is smaller than $r$, therefore the inductive hypothesis
gives the result.

\par {\sc Universit\'e de Grenoble I, Institut Fourier,
UMR 5582, B.P.74, 38402 Saint-Martin D'H\`eres Cedex, and IUFM de Lyon, 5 rue Anselme, 69317 Lyon Cedex (FRANCE)}\\ 

{\sc Department of Mathematics,
University of Ioannina, Ioannina 45110 (GREECE)}

\end{document}